%
%
\input amstex
\documentstyle{amsppt}

\expandafter\edef\csname aaaaa \endcsname{%
\catcode`\noexpand\@=\the\catcode`\@\space}
\catcode`\@=11
\def\nologo{\let\logo@\empty}
\csname aaaaa \endcsname
\nologo

\magnification=\magstep1
\parindent=1em
\baselineskip=18pt

\hsize=16 true cm
\vsize=24 true cm

\def\ove#1{\overline{#1}}
     \def\({\left(}       \def\al{\alpha}        \def\lee{\leqslant}
     \def\){\right)}      \def\e{\varepsilon}    \def\gee{\geqslant}
     \def\[{\left[}       
     \def\]{\right]}      
                          \def\be{\beta}
                                     
     \def\<{\langle}                 
     \def\>{\rangle}                 

\CenteredTagsOnSplits                        \NoBlackBoxes
\NoRunningHeads

    \topmatter
        \title {
      On Banach spaces without the approximation property$ {}^\dag$}
        \endtitle
          \author { O. I. Reinov}  \endauthor
\address\newline
  Dept. of Math. and Mech. of St.-Petersburg University \newline
  St.-Petersburg, St. Peterhof, Bibliotechnaja pl. 2 \newline
  198904, Russia, \newline
\endaddress
\email
orein\@orein.usr.pu.ru
\endemail
    \endtopmatter

\footnote" "{$ {}^\dag$
Reinov O.I., {\it On Banach spaces without the approximation property},\,
Funkts. analiz i ego prilozhen. (1982), v. 16, vyp. 4, p. 84-85
(in Russian).}
\document

{\bf 1.} If $ X$ is a Banach space of type $ p$ and of cotype $ q,$
then every its $ n$-dimensional subspace is $ Cn^{1/p-1/q}$-complemented
in $ X$ (cf. [1]). Szankowski [2] has showed that if
$ T(X)=\sup \{ p:\ X \text{ of type } p \}\neq2$ or
$ C(X)=\inf \{ q:\ X \text{ of cotype } q \}\neq2,$  then $ X$ has
a subspace without the approximation property. Thus, if each
subspace of $ X$ possesses the approximation property, then necessarily
$ T(X)=C(X)=2$ and, therefore, all of finite dimensional subspaces
in $ X$ are "well" complemented. Moreover, the space $ X$ need not be
Hilbertian (or isomorphic to a Hilbert space). The examples of such
spaces were constructed by Johnson [3].

In connection with the examples of Szankowski and Johnson,
the natural questions arises:
1){\it if $ T(X)=C(X)=2,$ then is it true that every subspace
in $ X$ has the approximation property}? \, 2) how "well complemented"
may be the finite dimensional subspaces of a space without
the approximation property: in particular, {\it does there exist
a space $ X$ without the approximation property, the constants $ C$ and $ A$
such that each $ n$-dimensional subspace of $ X$\, is
$ C\,\log^A n$-complemented}?
This note is devoted to the answers to these questions
(negative for the first one and positive for the second).
\medpagebreak

{\bf 2.} Our example is based entirely on the construction of
Szankowski [2].
Let $ I_n= \{ 2^n+1,\dots , 2^{n+1}\};$ $ \{ e_k\}_1^\infty$ be
the standard basis in $ c_0;$ $ \{ e'_k\}_1^\infty$ be a sequence of
the corresponding coordinate functionals. Further, let
$ z_i=e_{2i}-e_{2i+1}+e_{4i}+e_{4i+1}+e_{4i+2}+e_{4i+3}$ and
$ z_i'= 2^{-1}(e'_{2i}-e'_{2i+1}).$ Denote by $ W$ the linear span
of the set $ \{ z_i\}$ and put, for linear maps
$ U: W\to W$ and for $ n=1,2,\dots,$  \,
$\be_n(U)=2^{-n}\sum_{I_n} \< z'_i, Uz_i\>.$

In [2], there were constructed, for $ n=1,2,\dots,$ the partitions
$ \Delta_n$ and $ \nabla_n$ of the set $ I_n$ and finite scalar
sequences $ y_j\in c_0 \,(j=1,2,\dots)$ with the following properties:

a) $\be_n(U)-\be_{n-1}(U)= 2^{-n-1}\sum_{I_n} \< e'_i, y_i\>;$

b) if $ B_1, B_2\in\nabla_n,$ then $ \ove{\ove B}_1 =\ove{\ove B}_2=
       m_n\gee  C2^{n/8};$ hence, $ \ove{\ove\nabla}=2^n/m_n;$

c) if $ A\in\Delta_n, B\in\nabla_n,$ then $ \ove{\ove{B\cap A}}\lee 1,
         \ove{\ove A}\gee  C2^{7n/8};$

d) if $ B\in\nabla_n,$ then $ \sum_{i\in B} \e_iy_i=
     \sum_{j=1}^{10} \( \pm\sum_{i\in A_j} \e_ie_i\),$ where
$ \e_i=\pm1;$ $ A_j (j=1,\dots, 10)$ are subsets of some elements from
$ \Delta_{n-1},$ $ \Delta_{n}$ and $ \Delta_{n+1}$ and, moreover,
$ \ove{\ove {A_j}}=\ove{\ove B}.$

It follows from a) and  b) that

e) $ \be_n(U)-\be_{n-1}(U)=
  2^{-n-1}\sum_{B\in\nabla_n} 2^{-m_n} \sum_{\e}
\< \sum_{i\in B} \e_ie'_i, \sum_{i\in B} \e_iUy_i\>,
$
where $ \sum_\e$ denotes the summing over all collections of the signs
$\e_i=\pm1 $ $ (i\in B).$
\medpagebreak

{\bf 3.} An example of a space without the approximation property
which has good enough finite dimensional subspaces will be found
among the subspaces of the space
$ Y=\( \sum_n \( \sum_{A\in\Delta_n} l^2_{\ove{\ove A}}
             \)_{l^{p_n}}
    \)_{l^2}, $
where the numbers $ p_n$ are defined from the relations
$ 1/p_n-1/2=n^{-1}\log_a n^{1+\e}$
$ (\e>0, a= \root{8}\of{2}).$ It is clear that $ T(Y)=C(Y)=2.$ Let
$ X_\e$ be the closure in $ Y$ of the linear subspace $ W.$

\proclaim {\bf Theorem}\it
The space $ X_\e$ does not possess the approximation property.
There exists a constant $ C_\e>0$ such that
if $ E$ is an $ n$-dimensional subspace of $ X_\e,$ then
$1)$\, $ d(E, l^2_n)\lee C_\e \log^{1+\e} n$ and
$2)$\,
$E$\ is  $C_\e\log^{1+\e} n$-complemented in $ Y$ (and, hence, in $ X_\e$).
\endproclaim\rm
\medpagebreak

{\bf 4. \it Proof.} If $ B\in\nabla_n$ then, for $ \e_i=\pm1,$
we get from b), c) and d):
$$ \multline \text{f) }\,
\bigg\|\sum_{i\in B} \e_ie'_i\bigg\|_{Y^*}\lee
\( \sum_{A\in \Delta_n} \ove{\ove{(A\cap B)}}^{\,p'_n/2}\)^{1/p'_n}\lee \\ \lee
   \ove{\ove{ \{ A\in\Delta_n:\ A\cap B\neq\emptyset\}}}^{\,1/p'_n}\lee
    \ove{\ove {B}}^{\,1/p'_n}=m_n^{1/p'_n};
\endmultline
$$
$$ \multline \text{g) }  \,
\bigg\|\sum_{i\in B} \e_iy_i\bigg\|_{Y}\lee
     10\,\max_{1\lee j\lee 10} \bigg\| {\shave\sum_{i\in A_j}} \e_ie_i
                              \bigg\| \lee
   10\,\max_{1\lee j\lee 10} \( \ove{\ove {A_j}}\)^{1/2}\lee 10 m_n^{1/2}.
\endmultline
$$
Thus, it follows from e)
$$ \multline \text{h) }\,
 |\be_n(U)-\be_{n-1}(U)|\lee  \\ \lee
  2^{-n-1} \ove{\ove{\nabla_n}} m_n^{1/p'_n}
  \max \left\{
     \bigg\| U \bigg( \sum_{i\in B} \e_iy_i\bigg)\bigg\|:\
\e_i=\pm1,\,B\in\nabla_n
\right\}.
\endmultline
$$

Set $ F_n= \{ m_n^{-1/p_n} \sum_{i\in B} \e_iy_i:\
      \e_i=\pm1,\,B\in\nabla_n   \}.$
Then $ |\be_n(U)-\be_{n-1}(U)|\lee 2^{-1} \max
  \{ \|Uf\|:\ f\in F_n\}$ and, moreover,
$$ \multline
 \al_n=\max \{ \|f\|:\ f\in F_n\} \lee  \\ \lee
                m_n^{-1/p_n}\cdot 10 m_n^{1/2}\lee
   10 (C2^{n/8})^{-n^{-1}\,\log_a n^{1+\e}} \lee
  \operatorname{ const} \cdot n^{-1-\e}.
  \endmultline
$$

So we are in conditions of Lemma 1 of [2] and, therefore,
the space $ X_\e$ does not have the approximation property.

For $ T\in \L(X,Z),$ let us put $ \gamma_2(T)=\inf \{ \|A\|\,\|B\|:\
T=BA: X\to l^2\to Z\}.$ Let
$$ Y_N=
\( \sum_{n\lee N} \bigg( \sum_{A\in\Delta_n} l^2_{\ove{\ove A}}\bigg)_{l^{p_n}}
\)_{l^2}\ \text{ and }\
   Y^N=
\( \sum_{n> N} \bigg( \sum_{A\in\Delta_n} l^2_{\ove{\ove A}}\bigg)_{l^{p_n}}
\)_{l^2}.
$$

Let $ P_N$ and $ P^N$ be the natural projections from $ Y$ onto $ Y_N$ and
$ Y^N$ respectively. Since $ \ove{\ove{\Delta_n}}\lee C 2^{n/8},$ then
$$ d \( \bigg( \sum_{A\in \Delta_n} l^2_{\ove{\ove A}}\bigg)_{l^{p_n}},
      l^2_{2^n}\)\lee  (\ove{\ove{\Delta_n}})^{\,1/p_n-1/2}\lee C_1n^{1+\e};
$$
Whence, the Banach-Mazur distance from $ Y_N$ to an euclidean space
does not exceed $ C_2N^{1+\e}.$

Let now $ E$ be an arbitrary $ 2^N$-dimensional subspace of $ Y;$
$ E_N=P_{8N}(E)$ and $E^N=P^{8N}(E).$
Because of above arguments, there exists a projection
$ P_1$ from $ Y_{8N}$ onto $ E_N,$ and moreover
$ \gamma_2(P_1)\lee C_3N^{1+\e}.$
On the other hand, $ Y^{8N}$ is a space of cotype 2 and of type $ p_{8N};$
therefore, there exists a projection $ P_2$ from $ Y^{8N}$ onto $ E^N$
such that
$ \gamma_2(P_2)\lee C_4 (\dim E^N)^{1/p_{8N}-1/2}\lee C_5 N^{1+\e}.$
Put $ P_0= P_1P_{8N}+P_2P^{8N}$ and $ E_0=P_0(Y).$
Then $ P_0$ is a projection from $ Y$ onto $ E_0,$ and
$ \gamma_2(P_0)\lee C_6N^{1+\e}.$
Hence, firstly, there exists a projection from $ Y$ onto $ E$
with the norm $ \lee C_6 N^{1+\e}$
and, secondly, $ d(E, l^2_{2^N})\lee C_6 N^{1+\e}.$
\medpagebreak

{\bf 5.} If we define $ \{p_n\}$ from the relations
$ p_n^{-1}-2^{-1}=n^{-1}\, \log_an^{1+\e_n},$  where
$ \e_n\to 0$ and the series $ \sum n^{-1-\e_n}$ is convergent,
then the same arguments show that
{\it  there exists a Banach space $ X$ without the approximation property
and such that if $ \e>0$ and $ E$ is a $ n$-dimensional subspace in $ X,$
then there is a projection $ P$ from $ X$ onto $ E$
such that  $ \gamma_2(P)\lee C_\e\,\log^{1+\e} n.$ }
In particular, one can take, for each $ \delta>0,$ \,
$ n^{\e_n}= (\log n)(\log\log n)\dots (\log^{1+\delta} \log\dots\log n).$
Naturally, a question arises:
does there exist a Banach space without the approximation property,
all $ n$-dimensional subspaces of which are
$ C\log n$-complemented? $ C\log n$-euclidean?



\Refs

\ref \no 1  \by Pisier G.  \pages 1--21
\paper Estimations des distances \`a un espace euclidien et
   des constantes de proj\'ection des espace de Banach
   de dimensions finie
\yr 1979 \vol  \issue exp. 10,\nofrills
\jour  Seminaire d'analyse fonctionelle 1978--1979
\endref

\ref \no 2  \by A. Szankowski  \pages 123--130
\paper  Subspaces without approximation property
\yr 1978 \vol 30 \issue
\jour Israel J. Math.
\endref

\ref \no 3  \by Johnson W.B.  \pages  1--11
\paper Banach spaces all of whose subspaces have the approximation property
\yr 1980 \vol  \issue exp. 16,\nofrills
\jour Seminaire d'analyse fonctionelle 1979--1980
\endref

\endRefs
\enddocument